\documentclass[12pt,oneside]{amsart}
\usepackage{amsxtra}
\usepackage{amsopn}
\usepackage{amsmath,amsthm,amssymb}
\title{Hyper-K\"ahler quotients of solvable Lie groups}
\author{M. L. Barberis}
\author{ I. Dotti}
\author{ A. Fino}
\thanks{Research partially supported by
Conicet, Antorchas, SecytU.N.C (Argentina), by MIUR, CNR (Italy) and  ESI (Vienna)}
\address{Laura Barberis and Isabel G. Dotti: FaMAF, Universidad
Nacional de C\'ordoba\\5000 C\'ordoba,
Argentina}
\email{barberis@mate.uncor.edu, idotti@mate.uncor.edu}
   \address{Anna Fino: Dipartimento di Matematica, Universit\`a di Torino\\
Via Carlo Alberto 10, 10123 Torino, Italy}
\email{fino@dm.unito.it}

\newtheorem{teo}{Theorem}[section]
\newtheorem{prop}{Proposition}[section]
\newtheorem{corol}{Corollary}[section]

\theoremstyle{definition}

\theoremstyle{remark}
\newtheorem*{rem}{Remark}

\newcommand{\beq}{\begin{equation}}
\newcommand{\eeq}{\end{equation}}
\newcommand{\bqn}{\begin{eqnarray}}
\newcommand{\eqn}{\end{eqnarray}}
\newcommand{\bqne}{\begin{eqnarray*}}
\newcommand{\eqne}{\end{eqnarray*}}

\newcommand{\R}{{\Bbb R}}
\renewcommand{\H}{{\Bbb H}}

\newcommand{\xom}{\mbox{\boldmath$\Omega$}}

\begin{document}

\begin{abstract} In this paper we apply the   hyper-K\"ahler quotient
construction to Lie groups with a left invariant hyper-K\"ahler structure under the action of
a closed abelian subgroup  by left multiplication. This is motivated by the fact that 
some known hyper-K\"ahler metrics can be recovered in this way by considering 
different Lie group structures on $\H^p \times \H^q$ ($\H$: the quaternions). We obtain new complete hyper-K\"ahler metrics on
Euclidean spaces and give their local expressions.
\end{abstract}

\maketitle

\section{Introduction}

Hyper-K\"ahler reduction allows to construct hyper-K\"ahler manifolds from others that admit a group acting by tri-holomorphic
isometries \cite{HKLR}. Families of $4n$ dimensional hyper-K\"ahler quotients admitting a
tri-holomorphic $T^n$-action were constructed in \cite{G,BD}. In particular, in \cite{BD} the geometry and topology of
hyper-K\"ahler quotients of $\H^d$ by subtori of $T^d$ has been studied.

The hyper-K\"ahler quotient construction has been also applied in \cite{GRG} to the flat space $\H^d$ to obtain some monopole
moduli space metrics in explicit form using \cite{LR,PP}, for instance the Taubian-Calabi \cite{Ro}  and the  Lee-Weinberg-Yi
metric  \cite{LWY}. These are constructed  by considering the following actions of $\R$ on
$\H
\times
\H^m$ (resp.
$\R^m$ on  $
\H^m \times \H ^m$):
\begin{eqnarray*}
\R \times \H \times \H^m &\to &\;\; \H \times \H^m\\
(t, (q,w_1, \ldots, w_m)) &\to &(t+ q, e^{i t} w_1,
\ldots, e^{i  t} w_m), 
\end{eqnarray*}
 \begin{eqnarray*}
\R^m \times \H^m \times \H ^m &\to &\;\; \H^m \times \H ^m  \\
((t_1, \ldots, t_m), (q_1, \ldots, q_m, w_1, \ldots, w_m)) &\to &(t_1 + q_1
, \ldots,  t_m + q_m ,\\{}&{}& e^{i \langle \theta_1 ,T \rangle} w_1, \ldots, 
e^{i \langle \theta_m , T \rangle } w_m), 
\end{eqnarray*}
where $\theta \in \text{GL}(m,\R)$, 
$T = (t_1, \ldots, t_m), \; \theta_{\beta}$ are the rows of $\theta$ and 
$\langle \,\, , \, \rangle$ is the Euclidean inner product in $\R^m$. 
The first  action  gives rise to the Taubian-Calabi metric, which 
coincides with the Taub-Nut metric
for $q=1$, and the second one corresponds to  the Lee-Weinberg-Yi metric. We show that in 
both cases the  metric can be recovered by endowing $\H \times \H^m$  (resp. $\H^m \times \H^m$) 
 with a hyper-K\"ahler Lie group structure and taking the quotient with respect to
 a suitable closed abelian subgroup. 

In the present work we study hyper-K\"ahler quotients starting from a Lie group $G$ with a left invariant hyper-K\"ahler
structure.
Such a group  is necessarily flat since it is Ricci flat and
homogeneous (see \cite{AK}).   It follows from \cite{M} that  $G$
must be $2$-step solvable and when $G$ is simply connected, $G$ is a
semidirect product of the form
$\H^p \ltimes_{\theta} \H^q$, where $\theta$ is a homomorphism 
from $\H^p$ to $T^q$, a maximal torus in Sp$(q)$  (see Proposition \ref{hyperkahlerflat}
 and \eqref{homtheta}). This leads us to get a characterization of hyper-K\"ahler Lie groups. 

We  take a connected closed abelian subgroup
$\R^l$ ($l \leq p$) of  $\H^p$ which acts
  on $G$ by left translations, hence the action is free and the 
moment map has no critical points. 
 This
action is tri-Hamiltonian, therefore the hyper-K\"ahler quotient construction \cite{HKLR} can be applied.
We prove that the metric obtained on the hyper-K\"ahler
quotient is complete and the quotient is diffeomorphic to an Euclidean space. Since the $\R^l$-action commutes  with an action
of 
  the torus $T^q$, if $l = p$ the
$4q$-dimensional hyper-K\"ahler quotient admits a tri-holomorphic
$T^q$-action. Such action  has a unique fixed point when $p = q$. In this way we obtain new complete hyper-K\"ahler metrics
which generalize the Taubian-Calabi  and the Lee-Weinberg-Yi metrics.  Using the same method as in
\cite{GRG,PP,LR}, we obtain a local
expression  of the hyper-K\"ahler quotient metrics. This expression is given in terms of the structure constants of the
corresponding  Lie group $\H^p \ltimes_{\theta} \H^q$.

\smallskip

{\it Acknowledgments.} The authors thank Andrew Swann for  useful comments and Nigel Hitchin for helpful suggestions. The
first and third authors are grateful for the hospitality at the International Erwin Schr\"odinger Institute for
Mathematical Physics in Vienna during the program \lq \lq Geometric and analytic problems related to Cartan connections \rq
\rq.

\section{ Preliminaries}

Let $(\mathfrak g , g)$ be a metric Lie algebra, that is, $\frak g$ is a Lie algebra endowed with an  
inner product $g$.  The Levi Civita connection associated to the metric
can be computed by
\begin{eqnarray}\label{LC}
2g(\nabla _XY,Z) &=&g([X,Y],Z)
-g([Y,Z],X) + g([Z,X],Y)),
\end{eqnarray}
for any $X,Y,Z$ in $\mathfrak{g}$.

A hypercomplex structure on $\mathfrak g$ is a triple of complex structures
   $\{J_{\alpha}\}_{\alpha=1,2,3}$ satisfying the quaternion relations
$$J_{\alpha}^2= -\text{id},\; \alpha=1,2,3,\qquad  J_1J_2=-J_2J_1=J_3,$$
together with the vanishing of
   the Nijenhuis tensor
$N_{\alpha}(X,Y)= 0$,
for any $X,Y\in \mathfrak g$ and ${\alpha}=1,2,3$.  Here, the Nijenhuis tensor
stands for  
$$N_\alpha(X,Y)= J_\alpha([X,Y]-[J_\alpha X,J_\alpha Y])-([J_\alpha X,Y]+[X,J_\alpha Y]),$$
where $X,Y\in \mathfrak g$.

Let $\mathfrak g$ be a Lie algebra endowed with a hypercomplex
structure $\{J_{\alpha}\}_{{\alpha}=1,2,3}$ and an inner product $g$,
compatible with the hypercomplex structure, that is
$$
g (X,Y) = g (J_1X,J_1Y) = g (J_2X,J_2Y)= g(J_3X,J_3Y),
$$
for all $X,Y \in \mathfrak g.$   We will say that $(\mathfrak g , \{
J_{\alpha} \}, g)$ is a
hyper-K\"ahler Lie algebra when $(\mathfrak g , 
J_{\alpha} , g)$ is a
K\"ahler Lie algebra, for each $\alpha $, that is, 
$\nabla  J_{\alpha}=0$, where $\nabla$ is the Levi-Civita
connection of $g$.   
 This is equivalent to $d \omega_{\alpha} =0$, where 
$\omega_{\alpha}$ are the associated K\"ahler forms defined by
 $\omega_{\alpha}(X,Y)= g(J_\alpha X,Y), \; X,Y \in \mathfrak g.$

   If $G$ is a Lie group with Lie algebra $\mathfrak g$ then the above
structures
   on $\mathfrak g$ can be left  translated to all of $G$ obtaining
invariant hyper-K\"ahler structures on $G$.

Note that a Lie group with an invariant hyper-K\"ahler structure is
necessarily flat since a hyper-K\"ahler metric is Ricci flat and in
the homogeneous case, Ricci flatness implies flatness (see \cite{AK}).
Examples of non commutative Lie groups carrying a flat  invariant metric
  are given by $T^k \ltimes \R^m$ where $T^k$ is a torus in SO$(m)$. The next
  proposition, which is a consequence of the characterization of flat Lie
algebras given in \cite{M}, shows that this family of examples
essentially exhausts
  the class (see also \cite{DM}). This will allow us to
give a characterization of hyper-K\"ahler Lie algebras as a special
class of subalgebras of  $\R^s \times \frak{e}(4q)$, where
$\frak{e}(4q)
= \frak{so}(4q)\ltimes \R^{4q}$ is the Euclidean Lie algebra.

\begin{prop} [\cite{M}] \label{flat} Let $(\frak g , g)$ be a flat Lie algebra. Then
$\frak g$ decomposes orthogonally as $$\frak g=
\frak z (\frak g)\oplus \frak h \oplus \frak g ^1,$$ where ${\frak z}
({\frak g})$ is the center of ${\mathfrak g}$,  $\frak h$
is an abelian Lie subalgebra,
 the commutator ideal ${\frak g}^1$  is abelian and
  the following conditions are satisfied:

i) ad$: \frak h \to \frak{so}(\frak g ^1)$ is injective and ${\frak g}^1$ is even dimensional;

ii) ad$_X= \nabla _X $ for any $X \in \frak z (\frak g) \oplus {\mathfrak h}
$.

In particular, $\frak g$ is isomorphic to a Lie subalgebra of  $\R^s
\times \frak e (\frak g ^1)$, where $s=\dim \frak z (\frak g)$.

\end{prop}

\begin{proof} By \cite{M} a flat Lie algebra $(\mathfrak g , g)$ decomposes
orthogonally as
\begin{equation} \label{split}
{\mathfrak g} = {\mathfrak h} \oplus {\mathfrak b},
\end{equation}
where $\mathfrak h$ is an abelian Lie subalgebra, $\mathfrak b$ is the
abelian ideal defined by $\{B \in {\mathfrak g}: \nabla_B=0\}$ and
$${\mbox {ad}}_X: {\mathfrak b} \to {\mathfrak b}$$
is skew-symmetric, for any $X \in {\mathfrak h}$.  We observe that
the choice of $\frak b$, the fact
that $\nabla$ is torsion free
and that $\frak h$ is abelian imply
\begin{equation} \label{ad=nabla'}{\mbox {ad}}_X= \nabla_X, \qquad
\text{ for any } X \in {\mathfrak h}.\end{equation}
The above equation and the choice of $\frak b$ imply
\begin{equation}\label{injective}
{\mbox{ad}}: \frak h \to \frak{so}(\frak g )\end{equation} is injective.

We notice next that the decomposition \eqref{split} implies that
  ${\mathfrak
g}^1 \subseteq {\mathfrak b}$,  hence $\mathfrak b$ decomposes orthogonally
as
$$
{\mathfrak b} =  {\mathfrak v}\oplus {\mathfrak g}^1.
$$
   We show below
  that ${\mathfrak v} = {\mathfrak z}({\mathfrak g})$, where
${\mathfrak z}({\mathfrak g})$ denotes  the center of $\mathfrak g$.
In particular,
$\frak g$ will decompose orthogonally as $$\frak g=
\frak z (\frak g)\oplus \frak h \oplus \frak g ^1,$$ with  $\frak h$
and $\frak g ^1$  abelian and such that~$i)$ holds. To show that 
$\frak g^1$ is even dimensional, assume that $\dim \frak g^1 =2m+1$. 
Since ad$_X, \; X \in \frak h$, is a commutative family of endomorphisms
 in $\frak{so} (2m+1)$, they are conjugate to elements in a maximal abelian 
subalgebra of $\frak{so} (2m+1)$, hence there exists $Z \in \frak g^1$ 
such that ad$_X (Z) =0$ for any $X \in \frak h$, therefore $Z \in \frak z (\frak g )
\cap \frak g^1$, a contradiction. 

Since ${\mbox {ad}}_X: {\mathfrak g}^1 \to {\mathfrak g}^1$ and it
is skew-symmetric, for any $X \in {\mathfrak h}$,  it preserves
${\mathfrak v}$. Therefore, $[X, \frak v] \subset \frak v \cap \frak
g ^1
= 0$  for $X\in \frak h$ and $\frak v \subset  {\mathfrak z}({\mathfrak g})$
follows. On the other hand, if $Y\in {\mathfrak z}({\mathfrak g})$, then:
\[   0= g([Y,X],U)= g(Y, [X,U]),  \]
for every  $X\in \frak h, \; U \in \frak g ^1$, that is,
${\mathfrak z}({\mathfrak g}) \perp \frak g ^1$ since $\frak g ^1 =
[\frak h , \frak g ^1]$. From  ${\mathfrak z}({\mathfrak g})
\cap \frak h = 0$ one has that $\frak v =
{\mathfrak z}({\mathfrak g})$.

Finally, using \eqref{LC} one can compute $\nabla_Y=0$ for $Y \in
{\mathfrak z}({\mathfrak g})$.  This together with \eqref{ad=nabla'}
imply~$ii)$ and the proposition follows.

\end{proof}

We will say that two  flat Lie algebras $(\frak g_1 , g_1)$ and $(\frak g_2 , g_2)$
 are  equivalent if there exists an orthogonal Lie algebra isomorphism $\eta : \frak g _1
\to \frak g_2$. Note that $\eta: \frak z (\frak g _1 ) \to \frak z (\frak g _2 )$, 
$\eta : \frak g _1^1 \to \frak g _2 ^1$ and therefore $\eta: \frak h_1 \to \frak h_2$ 
(see Proposition~\ref{flat}). Let ad$_i: \frak h_i \to \frak{so}(\frak g ^1_i)$, $i=1,2$,  be the 
corresponding injective maps induced by the adjoint representation on $\frak g_i$. 
Then the following diagram is commutative:
\[ \begin{array}{ccc}
\frak h _1 &\stackrel{\text{ad}_1}{\longrightarrow} &  \frak{so}(\frak g ^1_1) \\
\eta \Big\downarrow &   & \Big\downarrow  I_{\eta}\\
\frak h _2 &\stackrel{\text{ad}_2}{\longrightarrow} &  \frak{so}(\frak g ^1_2)
\end{array}
\]
where $I_{\eta}$ denotes conjugation by $\eta$. Conversely, 
it follows from Proposition \ref{flat} that every flat Lie algebra  with $2m$ dimensional commutator and $s$ dimensional 
center is equivalent to $ \R^s \times \R^k \ltimes _{\rho} \R^{2m}$, where 
  $\rho : \R ^k \to \frak{so}(2m)$ is injective, $\rho (\R^k)\R^{2m}=\R^{2m}$, the 
only non zero Lie brackets being $$[X, Y]=\rho(X)Y,\; X \in \R^k , \; Y \in \R^{2m}.$$

Given a flat Lie algebra $\R^s \times (\R^k \ltimes _{\rho} \R^{2m})$, the family 
  $\{\rho(T):  T \in \R^k \} \subseteq \frak{so}(2m) $ 
is an abelian subalgebra, then it is conjugate by an element in $SO(2m)$ to a subalgebra of the following 
maximal abelian subalgebra of $\frak{so}(2m) $:
\[  \frak t ^m= \left\{  \begin{pmatrix} 0 & -\phi_1 & & &  \\
\phi_1 &0 & & & \\
 &   &   \ddots &  & \\
  &  &  & 0& -\phi_m \\
  &  &  & \phi_m & 0
\end{pmatrix} : \phi_{\alpha} \in \R  \right\}           \]
with respect to an orthonormal basis $\{f_1 , \ldots , f_{2m}\}$   of $\R^{2m}$. 
In particular, $k \leq m$ and we may assume  
that 
any flat Lie algebra is equivalent to a Lie algebra 
such that  $\rho (\R ^k) \subset \frak t ^m$.   

Let   $\theta =(\theta _{\beta}^{\alpha})$ be the real $m \times k$  matrix of rank $k$ such that 
\begin{equation} \label{eqmu}  \rho(e_{\alpha})= \begin{pmatrix} 0 & -\theta_1^{\alpha} & & &  \\
\theta_1^{\alpha} &0 & & & \\
 &   &   \ddots &  & \\
  &  &  & 0& -\theta_m^{\alpha} \\
  &  &  & \theta_m^{\alpha} & 0
\end{pmatrix}, \quad 1\leq \alpha \leq k,  \end{equation}
 where
$\{ e_1, \ldots , e_k \}$ is an orthonormal 
basis of $\R^k$. The  condition $\rho (\R^k)\R^{2m} =\R^{2m}$ is equivalent to the fact that 
every row $\theta _{\beta}$ of $\theta$ is non zero. 

We introduce some notation that will be used in the next result. 
Let $M(k,m;k)$ be the set of $m \times k$ real matrices of rank $k$. $M(k,m;k)$
can be viewed inside  End$(\R^k,\frak{so}(2m))$   by means of the inclusion $\rho$:
\[    M(k,m;k) \hookrightarrow    \text{End}(\R^k,\frak{so}(2m)), \qquad \theta \mapsto \rho _{\theta}. \]
We identify  $M(k,m;k)$ with its image under $\rho$ and let  
$O(k)\times O(2m)$ act on $M(k,m;k)$ as follows:
\begin{equation} \label{flatequiv} O(k)\times O(2m)  \times M(k,m;k) \to M(k,m;k), \qquad  
(A,B, \; \rho_{\theta}) \mapsto B \rho _{(\theta A)}B^{-1},  \end{equation}
where $ B \rho _{\theta }B^{-1} \in  \text{End}(\R^k,\frak{so}(2m))$ is defined 
by $ B \rho _{\theta }B^{-1} (T)=  B \, \rho _{\theta }(T)\, B^{-1}$, $T \in \R^k$.  
It follows from the definition of equivalence between flat Lie algebras 
that 
$$\R^k \ltimes _{\rho_{\theta}}\R^{2m} \cong \R^k \ltimes _{\rho_{\theta '}}\R^{2m} $$ 
if and only if $\rho _{\theta}$ and $\rho_{\theta '}$ lie in the same $O(k)\times O(2m)$-orbit. %In other words, 
%the space of flat Lie algebras $(\frak g , g)$ with trivial center such that $\dim g= k+2m$, %$\dim \frak g ^1 \leq 2m$ is parameterized by the the orbit space  $M(k,m;k)/ O(k)\times %O(2m)$. 

The next proposition summarizes the above results and gives the classification  of flat Lie algebras that will be needed 
in the next section (see also \cite{KTV}). 
\begin{prop} \label{matrixflat} Let $(\frak g , g)$ be a flat Lie algebra, $\dim \frak g ^1 =2m$, 
$ \dim \frak z (\frak g )=s$. Then  there exists  $\theta =(\theta _{\beta}^{\alpha})\in M(k,m;k)$  such that $\theta _{\beta} \neq 0$ for every $1\leq 
\beta \leq m$ and   
 $\frak g$ decomposes orthogonally as 
\[ \frak g \cong \R^s \times  (\R^k \ltimes _{\rho _{\theta }} \R^{2m}), \] where $\R^k \ltimes _{\rho_{\theta }} \R^{2m} $ has
an orthonormal basis $\{e_1, \ldots , e_k, f_1, \ldots , f_{2m}\}$  and $T \in \R^k $ acts on $\R^{2m}$ in the following way:
 \begin{equation} \label{eqflat}   \rho _{\theta}(T)= \begin{pmatrix} 0 & -
\langle T, \theta _1 \rangle &   &   &    \\
\langle T , \theta _1 \rangle & 0 & & &  \\
  &   &   \ddots & &   \\
  & &  &  0 & -\langle T , \theta _m \rangle  \\
 &  &  & \langle T , \theta _m \rangle & 0 
\end{pmatrix}  
\end{equation} 
where  $\langle \,\, , \, \rangle$  
denotes the Euclidean inner 
product on $\R^k$. Moreover, 
 \[ \R^k \ltimes _{{\rho}_{\theta}} \R^{2m} \cong \R^k \ltimes _{\rho _{\theta '} } \R^{2m}  \]
 as flat Lie algebras if and only if $\rho_{\theta}$ and $\rho _{\theta '}$ lie in the same 
$ O(k)\times O(2m)$-orbit under the action \eqref{flatequiv}.
\end{prop}

\begin{rem}
Note that the Lie algebra $\R^k \ltimes _{{\rho}_{\theta}} \R^{2m}$ is a Lie subalgebra 
of the Euclidean Lie algebra $\frak e (2m)$:
\[ \R^k \ltimes _{{\rho}_{\theta}} \R^{2m} \hookrightarrow \frak e (2m), 
\quad (T,W) \mapsto 
\begin{pmatrix} \rho _{\theta} (T) & W \\
0 & 0 \end{pmatrix},
\]
$T \in \R^k, \; W \in \R^{2m}$. However, the inner product on $\R^k \ltimes _{{\rho}_{\theta}} \R^{2m}$ is not the one induced
from $\frak e (2m)$. 
\end{rem}
The next corollary follows from the description given in Proposition \ref{matrixflat}. 

\begin{corol} Any even dimensional flat Lie algebra is K\"ahler flat. 
\end{corol}

\begin{proof}
Let 
$\frak g _{\theta}= \R^s \times (\R^k \ltimes _{{\rho}_{\theta}} \R^{2m})$ be as in Proposition \ref{matrixflat} and $J$  the
 endomorphism of $\frak g_{\theta}$ such that $J^2=-$id, $J f_{2i+1} =f_{2i}$, $i=0, \ldots , m-1$, and 
$J$  is orthogonal preserving $\R^s \times \R^k$. The integrability of $J$, that is, the vanishing of $N_J$, follows from 
$\rho _{\theta} (T) J = J \rho _{\theta} (T)$, for any $ T \in \R^k$. Moreover, 
$\nabla J=0$ since $\nabla _T = \rho _{\theta} (T)$, for $T \in \R^k$. Therefore 
$(\frak g ,  J, g)$ is K\"ahler flat. 
 \end{proof}

\section{Hyper-K\"ahler Lie groups}

We apply Proposition~\ref{flat} to give a characterization of the Lie algebras 
carrying a hyper-K\"ahler structure $(\{J_{\alpha} \}, g)$. 

\begin{prop} \label{hyperkahlerflat}
Let $({\mathfrak g}, \{J_{\alpha} \}, g)$, $\alpha
=1,2,3$, be a hyper-K\"ahler
Lie algebra.
Then $\mathfrak g$ decomposes orthogonally as
$$
{\mathfrak g} = {\mathfrak t} \oplus {\mathfrak g}^1, \qquad \frak z (\frak g) 
\subset \frak t ,
$$
with both ${\mathfrak t}$  and  ${\mathfrak g}^1 $ abelian and
$J_{\alpha}$-invariant, $\alpha=1,2,3$,  such that
\newline
i) ${\mbox {ad}}_X J_{\alpha} = J_{\alpha} {\mbox {ad}}_X$, for any
$X \in {\mathfrak t}$, $\alpha =1,2,3$;
\newline
ii) $g({\mbox {ad}}_X Y, Z) + g(Y, {\mbox {ad}}_X Z) = 0$, for any $X
\in {\mathfrak t}, Y, Z \in {\mathfrak g}$.
\end{prop}

\begin{proof} Since a hyper-K\"ahler Lie algebra is flat \cite{AK}, $\frak g$
decomposes orthogonally as \newline
$\frak g =\frak z (\frak g)\oplus \frak h
\oplus \frak g ^1$
  and the conditions of Proposition~\ref{flat} are satisfied.
   Set $${\mathfrak t} = \frak z (\frak g) \oplus {\mathfrak h}.$$
  We
show next that if $({\mathfrak g}, \{ J_{\alpha}\} , g)$, $\alpha
=1,2,3$, is a hyper-K\"ahler
Lie algebra then
${\mathfrak t}$ and ${\mathfrak g}^1$ are
$J_{\alpha}$-invariant, $\alpha =1,2,3$, and that condition $i)$ is satisfied.
Observe that if $ X\in {\mathfrak t}$ and $B\in {\mathfrak g}^1$,
using that $\nabla J_{\alpha} =0$ and $ii)$ of Proposition
\ref{flat}, one has
$$ J_{\alpha} [X, B]= J_{\alpha}\nabla_XB = [X, J_{\alpha}B],$$
therefore, $i)$ follows. Since $\frak g^1=[\frak h , \frak g ^1]$,
the above equation
  also implies  that $ \frak g ^1 $ is $J_{\alpha}$-invariant and the
decomposition $\frak t \oplus
  \frak g ^1$ satisfies the desired properties.

\end{proof}
We will say that two hyper-K\"ahler Lie algebras $(\frak g , \{ J_{\alpha} \}, g)$ and 
$(\frak g ' , \{ J_{\alpha}' \}, g')$ are equivalent if there exists 
an equivalence $\eta$ of metric Lie algebras such that $\eta J_{\alpha}=J_{\alpha}' \eta$, 
$\alpha =1,2,3$.

Consider the hypercomplex
structure on 
\[ \H^q=\{ (W_1, \dots , W_q) : W_{\alpha}=u_{\alpha}+y_{\alpha}i + z_{\alpha}j +w_{\alpha}k  \; :\; u_{\alpha}, y_{\alpha}, z_{\alpha}, w_{\alpha} \in \R  \} \]
 given by right multiplication by $-i, \, -j ,\, -k$:
$$
J_1= R_{-i}, \qquad \qquad J_2= R_{-j}, \qquad \qquad J_3=R_{-k}. 
$$
We identify $\H^q \cong \R^{4q}$ with the Euclidean metric and let Sp$(q)=O(4q) \cap$ GL $(q, \H)$, where  
\[ \text{GL} (q, \H) =\{ T\in  \text{GL}(4q, \R) : T J_{\alpha}=J_{\alpha} T,\; \alpha =1,2,3  \}. \] 
Let $\frak t ^q$ be the following maximal abelian subalgebra
of the Lie algebra ${\mathfrak {sp}} (q)$ of Sp$(q)$: \begin{equation} \label{torus} \frak t^q = 
\left \{ \; \begin{pmatrix}
  0&-\phi_1 &0&0&  &   & & &  \\
\phi_1 &0&0&0& & &  & & \\
0&0&0&-\phi_1 & & &  & &  \\
0&0& \phi_1 &0& &  &  & & \\
& & & &\ddots & & &  &  \\
   & &  & & &0&-\phi_q &0&0\\
& & &  & &\phi_q &0&0&0\\
& & &  & &0&0&0&-\phi_q \\
& & &  & &0&0&\phi_q &0 \end{pmatrix} \; : \; \phi_i \in \R \right \}. \end{equation} 
We obtain the analogue of Proposition \ref{matrixflat} using the same procedure 
as before. Observe that, in this case, $\R^s \times (\R^k \ltimes_{\rho_{\theta }} \H^q) \cong \R^s \times (\R^k \ltimes_{\rho
_{\theta '}} \H^q)$ as
 hyper-K\"ahler Lie algebras if and only if $\rho _{\theta }$ and 
$\rho _{\theta '}$ lie in the same Sp$(p) \times \text{Sp}(q)$-orbit, where $s+k=4p$ and the
action of Sp$(p) \times \text{Sp}(q)$ is the analogue of \eqref{flatequiv}.

\begin{prop} \label{hkflat} Let $(\frak g , \{ J_{\alpha} \}, g)$ be a hyper-K\"ahler Lie algebra 
with $\dim \frak g^1 =4q$ and $\dim \frak z (\frak g )=s$. Then there exists  $\theta =(\theta _{\beta}^{\alpha}) \in M(k,q;k)$ , $s+k =4p$,  such that 
$\theta_{\beta} \neq 0$ for $1\leq \beta \leq q$ and \[ 
\frak g \cong \R^s \times (\R^k \ltimes _{\rho _{\theta}} \H^q).    \]
 $\R^k \ltimes _{\rho _{\theta}} \H^q $ is the  Lie algebra with orthonormal basis
\[ \{e_1, \ldots , e_k, f_1, f_1 i , f_1 j, f_1 k, \ldots ,  f_q, f_q i , f_q j, f_q k \}  \]
 such that an element $T \in \R^k$ acts on $\H^q$ by  
\begin{equation} \label{HKLiealg}
   \rho_{\theta}(T)= 
\begin{pmatrix}\rho^1_{\theta}(T) & & \\
 & \ddots & \\
 & &  \rho^q_{\theta}(T) \end{pmatrix}, 
\end{equation}
where $\rho^{\beta}_{\theta}(T)$ is the following $4\times 4$ real matrix: 
 \[  \rho^{\beta}_{\theta}(T)= \begin{pmatrix} 0&- \langle T , \theta _{\beta} \rangle & 0& 0   \\
\langle T , \theta _{\beta} \rangle &0&0 & 0\\
 0&0 &0&-\langle T , \theta _{\beta} \rangle   \\
 0& 0& \langle T , \theta _{\beta} \rangle &0 
\end{pmatrix}  
\] and $\langle \, \, , \, \rangle$ denotes the Euclidean inner product on $\R^k$. 
The Lie algebra $\R ^s \times (\R ^k \ltimes _{\rho _{\theta}} \H^q)$ is hyper-K\"ahler with its natural hypercomplex
structure obtained by extending $R_{-i}, \; R_{-j}, \;R_{-k}$  on $\H^q$ with any pair of anticommuting complex endomorphisms
on $\R ^s \times \R ^k $  and the canonical inner product. Moreover, 
\[ \R^s \times (\R^k \ltimes_{\rho _{\theta}} \H^q) \cong \R^s \times (\R^k \ltimes_{\rho _{\theta '}} \H^q) \]
as hyper-K\"ahler Lie algebras 
 if and only if $\rho _{\theta }$ and 
$\rho _{\theta '}$ lie in the same Sp$(p) \times \text{Sp}(q)$-orbit.  
\end{prop}

\subsection{Examples}
 As a consequence of Proposition~\ref{hkflat} we have that 
 there is a one parameter family of $8$-dimensional
hyper-K\"ahler   Lie algebras $\frak g _{\theta}$:
\begin{equation}  \label{8dim} \frak g _{\theta } \cong \R^3 \times (\R\ltimes _{\theta} \H ),
\end{equation} where $\R\ltimes _{\theta} \H $ has an orthonormal basis $\{e_1, f_1, f_1 i , 
f_1j , f_1 k   \}$ and $e_1$ acts on $\H$ as follows:
$$\rho_{\theta}(e_1)=  \begin{pmatrix} 0 & -\theta & &  \\
 \theta &0 & & \\
  & & 0 & -\theta \\
   & & \theta & 0 
\end{pmatrix} .$$
Note that these are pairwise non equivalent  flat metric Lie algebras for different values 
of $\theta$, but they are isomorphic  Lie algebras for 
$\theta \neq 0$.

 In dimension $12$ 
there are infinitely many non isomorphic  Lie algebra structures admitting hyper-K\"ahler
 metrics. 
In fact, for a fixed real number $s\neq 0$ we define $\frak g _s =\R^3\times (\R \ltimes _s \H^2)$, where $\R \ltimes _s \H^2$
has an orthonormal basis as in the statement  of Proposition~\ref{hkflat} with $e_1 $ acting on $\H^2$ as follows:
 \[
\rho _s(e_1) = 
 \begin{pmatrix}
  0&-1 & & &  &   & &  \\
1 &0& & & & &  & \\
 & &0&-1 & & &  &  \\
 & & 1 &0& &  &  & \\
    & &  &  &0&-s & &  \\
& & &   &s &0& & \\
& & &   & & &0&-s  \\
& & &   & & &s &0 \end{pmatrix} .  \]

It turns out that $ \frak g _s$ and $ \frak g_{r} $
are non isomorphic for $s \neq r$. 

\medskip
 
We describe below the Lie bracket on $\frak g _{\theta} = \R^k \ltimes _{\rho_{\theta}} \H^q$: 
 \begin{equation}\label{bracket} \begin{split}[(X,W),(X',W')]=& (0, i(\langle X , \theta_1\rangle W'_1-\langle X' , \theta_1\rangle W_1), \dots , i(\langle X , \theta_q \rangle W'_q-\langle X' , \theta_q \rangle W_q))\\
=& (0,  \rho_{\theta}(X) W' -\rho_{\theta} (X') W),
\end{split} \end{equation} $X, X' \in \R ^k$.  

The product on the  simply connected Lie group $G_{\theta}=\R^k \ltimes _{ \theta } \H^{q}$ with
Lie algebra $\frak g_{\theta}$ is  given as follows:
\begin{equation} \label{product}
(X, W)\cdot (X',W')= (X+X' , \, W + \theta (X)  W' ),
\end{equation}
where
$ X, X' \in \R^k, \; W, W' \in \H^q$,  $\,
W'=(W'_1, \dots , W'_q)$ and 
\begin{equation}\label{homtheta} {\theta}(X) W'= (e^{i\langle X, \theta_1 \rangle} W'_1, \dots ,
e^{i\langle X, \theta_q \rangle} W'_q). \end{equation}
Using that $(X,W)^{-1}=-(X,\theta (-X)  W)$, conjugation by
 $(X,W)$ is given as follows:
\[  I_{(X,W)}(X',W') =  (X,W)\cdot (X',W')\cdot (X,W)^{-1}=(X', 
W+\theta (X ) W' - \theta(X')  W)    \]
and therefore
\begin{equation} \label{Adjoint} \text{Ad}(X,W)(X',W')= (X', \theta (X ) W') + [(0,W), (X', 0)] 
= (X', \theta (X ) W' - \rho _{\theta}(X')W),   \end{equation}
for $X, X' \in \R^k , \, W, W' \in \H^q$.
% where $\overline{W}=
%(\overline{W_1}, \dots , \overline{W_q})$ and $\overline{W_l}$ denotes 
%conjugation in $\H$. 

The left invariant flat metric $g$ on $\R^k \ltimes _{\rho _{\theta}} \H^{q}$ in coordinates 
$(x_1,\dots ,x_{k}, W_1,\dots ,W_q)$, where $W_j=(u_{j}, y_{j}, z_{j} , w_j)$, is
 the Euclidean metric   $$g=\sum_{j=1}^{k} dx_j ^2 + 
\sum_{j=1}^{q}(du_j^2 + dy_j^2 + dz_j^2+ dw_j^2).$$ 

We will need to express the Euclidean metric on $\H^q$ in coordinates adapted to the hyper-K\"ahler moment map. Any quaternion
may be written as
$$
W_{\beta} = e^{i \psi_{\beta}/2} a_{\beta}, \quad \beta = 1, \ldots, q,
$$
with $\psi_{\beta} \in (0, 4 \pi]$ and $a_{\beta}$ is pure imaginary, so that $\overline a_{\beta} = - a_{\beta}$.  Let 
$$
{\bf r}_{\beta} = \overline W_{\beta} i W_{\beta} = \overline a_{\beta} i a_{\beta} = - a_{\beta} i a_{\beta}.
$$
The flat metric on $\H^q$ in coordinates $(\psi_{\beta}, {\bf r}_{\beta})$, $\beta = 1, \ldots, q$, is given by
\begin{equation} \label{Gibbonsflat}
\frac{1} {4}  \sum_{\beta = 1}^q \left( \frac{1} {r_{\beta}} d {\bf r}_{\beta}^2 + r_{\beta} (d \psi_{\beta} + {\xom}_{\beta}
\cdot d {\bf r}_{\beta})^2 \right),
\end{equation}
where $$r_{\beta} = \vert {\bf r}_{\beta} \vert, \quad  \mbox{curl} (\xom_{\beta}) = \mbox{grad} \left( \frac{1}
{r_{\beta}} \right)$$ (the curl and grad operations are taken with respect to the Euclidean metric on $\R^3$ with cartesian
coordinates ${\bf r}_{\beta}$).

\section{Completeness of the hyper-K\"ahler quotient metric}

According to Proposition~\ref{hkflat} any simply connected 
 Lie group with a left  invariant hyper-K\"ahler 
structure is of the form  
$G _{\theta}= \R^s \times (\R^k \ltimes_{\theta} \H^q)$ ( $k\leq q$, $s+k=4p$) with the hyper-K\"ahler metric
$g= g_1 \times g_2$, where $g_1$ is the  Euclidean metric on 
 $\R^s \times \R^k$ and $g_2$ is the Euclidean metric on $\H^q$.
  Let $\frak g _{\theta}$ be the Lie algebra of $G _{\theta}$. The associated K\"ahler 
forms:
\[  \omega _{\alpha}((X_1, W_1), (X_2, W_2))= g(J_{\alpha} (X_1, W_1) , (X_2, W_2) ), \qquad (X_1, W_1), (X_2, W_2) \in \frak g _{\theta} ,\]
$\alpha =1,2,3$, when left translated to $G_{\theta} $ become:
\[ \omega _{\alpha}= \omega _{\alpha}^1 +
  \omega _{\alpha}^2, \]
where $\omega _{\alpha}^j$, $j=1,2, \; \alpha =1,2,3$ are the standard symplectic forms on a vector space. Therefore, $(G _{\theta}, g, \omega _{\alpha})$ is equivalent, as 
a hyper-K\"ahler
 manifold, to the product  \begin{equation} \label{hkprod}(\R^s \times \R^k , g_1, \{ \omega _{\alpha}^1 \}) 
\times (\H^q , g_2 , \{ \omega _{\alpha}^2 \} ).  \end{equation}

 We will apply the hyper-K\"ahler quotient construction 
due to \cite{HKLR} to the case when 
 $L $ is the connected closed abelian Lie subgroup $\R^l \subset 
\R^k$ with Lie algebra $\frak l =\text{span}_{\R} \{ e_1 ,\ldots , e_l \}$ such that 
$\frak l$ is isotropic with respect to $\omega _{\alpha}$ for each $\alpha$. The 
action of $L$   on $G_{\theta}$ will be given by left 
 translations, therefore it preserves the hyper-K\"ahler structure.   We recall next the quotient construction in our 
particular case. 

Let ${\mathcal X}_V$ be the vector field generated by the action of $L$,  
 that is,
the right invariant vector field such that ${{\mathcal X}_V} _e=V$, 
where $V \in
\frak l$.
  Observe  that  \[   0= L_{{\mathcal X}_V} \omega _{\alpha} = d(i({\mathcal X}_V)\omega _{\alpha}) +i({\mathcal X}_V) d \omega
_{\alpha},
\]
where $i({\mathcal X}_V)\omega _{\alpha}$ denotes the $1$-form obtained by taking the
interior product with ${\mathcal X}_V$. Since the action is symplectic with respect to
$\omega_{\alpha}$, $\alpha =1,2,3$, we have that $i({\mathcal X}_V)\omega _{\alpha}$, $\alpha = 1,2, 3$, is closed. 
$G_{\theta}$ is simply connected, thus  $H_{dR}^1(G_{\theta} ,\R)=\{ 0\}$ and  $i({\mathcal X}_V)\omega _{\alpha}$ is
exact, that is,
\[      i({\mathcal X}_V)\omega _{\alpha}= d \, {\mu^{\theta}_{\alpha}}^V,              \]
   where ${\mu ^{\theta} _{\alpha}}^V$ is a Hamiltonian function associated to $V$.
Putting all these functions together, we obtain a  map to the dual space of the Lie algebra of $L$
\[    {\mu ^{\theta}}_{\alpha}: G _{\theta}\to \frak l^*  \]
 defined by
\[     {\mu ^{\theta}} _{\alpha}(X,W)( V)= {\mu ^{\theta}} _{\alpha} ^V (X,W).
\]
$L$ acts on ${\mathfrak l}^*$ by the coadjoint action. When the ambiguities in the choices of ${\mu ^{\theta} _{\alpha}} ^V$
may be adjusted to make 
$\mu ^{\theta}_{\alpha}$ $L$-equivariant, one has
the hyper-K\"ahler moment map 
$$
\mu ^{\theta}: G_{\theta} \to {\mathfrak l}^* \otimes \text{Im} \, \H , $$
defined by $\mu ^{\theta}= \mu ^{\theta}_1 i + \mu ^{\theta}_2 j + \mu^{\theta}_3 k$. Our choice of $L$ implies that 
 $\mu ^{\theta}_{\alpha }$ is $L$-equivariant for each $\alpha$.  Indeed, 
the action $A$ of $L$ on $G_{\theta}$   given by left translations:
\begin{equation}\label{action} \begin{array}{ccl}  A: L \times
G_{\theta} & \to & \;\; G_{\theta} \\
   ((V,0) , (X,W))    & \to & (V,0) \cdot   (X,W) = (V+X , \, \theta (V)
W) \end{array}    \end{equation}
(recall \eqref{product})  can be viewed as a diagonal action of $L$:
$$A (V)(X,W)= (A_1(V)X , A_2 (V)W),$$
   where $A_1$ acts
by left translations on $\R^s \times \R^k$ and $A_2$ is a linear symplectic
action  on $\H^q$.
The  moment map $\mu ^{\theta}_{\alpha}$ corresponding to $A$ can be
obtained by adding up the moment
maps of $A_1$ and $A_2$ since \eqref{hkprod} holds. 
By a direct calculation one has (see \cite{GS}):
\[  \mu ^{\theta}_{\alpha}(X,W)( V) = \omega _{\alpha}(V,X) +
\frac 12 \omega _{\alpha}
(\rho _{\theta}(V) W, W)  \]
(see \eqref{bracket}).
The $L$-equivariance of the first term follows since $L$ is isotropic and the 
second term is $L$-equivariant since  it is the moment map of a linear action 
on a symplectic vector space (see \cite{GS}). 

Let $\xi \in \frak z \otimes 
\R^3$ be a regular value for $\mu ^{\theta}$, where $\frak z$ is the subspace of $\frak l ^*$ of invariant elements under the
coadjoint action,  and consider the quotient $L\backslash ({\mu^{\theta}} )^{-1}(\xi)$. 
According to \cite{HKLR}, when the action of $L$ on $({\mu ^{\theta}}) ^{-1}(\xi)$ is free with Hausdorff quotient manifold 
$L\backslash {(\mu ^{\theta}})^{-1}(\xi)$, 
the hyper-K\"ahler metric on $G_{\theta}$ induces a hyper-K\"ahler metric on $L\backslash 
({\mu ^{\theta}})^{-1}(\xi)$. 

Our hypotheses imply:
\begin{enumerate}
\item the center of ${\mathfrak l}^* $, that is, the  invariant elements under the coadjoint action, coincides with
${\mathfrak l}^* $ since $L$ is abelian;
\item  the action of $L$ on $G_{\theta}$ is free, hence it is free on  
$({\mu ^{\theta}})^{-1} (\xi)$, for any $\xi=\xi_1 i + \xi_2 j + \xi_3 k \in \text{Im}\, \mu ^{\theta}$. 
 In particular, any  $\xi \in \text{Im}\, \mu ^{\theta}$ is a regular 
value of  the hyper-K\"ahler moment map;
\item since $L$ is closed in $G_{\theta}$ and acts by left   translations, the set of  right cosets $L\backslash G_{\theta}$
is  a Hausdorff manifold and so is   
$L \backslash ({\mu ^{\theta}})^{-1}(\xi)$,  
  for any $\xi \in \text{Im}\, \mu ^{\theta} $. \end{enumerate}
Therefore, $L \backslash ({\mu ^{\theta}})^{-1}(\xi)$ inherits 
a hyper-K\"ahler metric.

We show next that the hyper-K\"ahler metric on $L \backslash ({\mu ^{\theta}})^{-1}(\xi)$ is 
always complete.  
The left  invariant
metric $g$ on $G_{\theta}$ induces in a natural way a metric $\tilde g$ on
$L \backslash G_{\theta}$ such that the natural projection
$$
\pi: (G_{\theta}, g) \to (L \backslash G_{\theta}, \tilde g)
$$
is a Riemannian submersion. The completeness of $g$
  implies that $\tilde g$ is also complete (see \cite{Her}).
Since $\xi$ is a regular value of $\mu ^{\theta}$, $({\mu ^{\theta}})^{-1} (\xi)$ is a closed
embedded submanifold of $G_{\theta}$ and one has the following commutative diagram:
$$
\begin{array} {ccc}
G _{\theta}& \stackrel{\pi}{\longrightarrow}&L \backslash G_{\theta}\\
\Big\uparrow&
&\Big\uparrow\\
{\mu ^{\theta}}^{-1} (\xi) &\stackrel{\tilde \pi}{\longrightarrow} &L \backslash
{\mu ^{\theta}}^{-1} (\xi),
\end{array}
$$
where the vertical arrows are the natural inclusions.

Since $\mu ^{\theta}$ is $L$-equivariant, it induces a map
$$
\tilde \mu^{\theta}: L \backslash G_{\theta} \to  \frak l^* \otimes \text{Im} \H
$$
and $ L \backslash {\mu ^{\theta}}^{-1} (\xi) = \tilde {\mu ^{\theta}}^{-1} (\xi)$. Then $
L \backslash {\mu ^{\theta}}^{-1} (\xi)$ is a closed $d$-dimensional embedded submanifold of
$L \backslash G_{\theta}$ ($d = \dim G_{\theta} - 4 \dim L$) and the
induced metric is complete. 
By O'Neill formula (\cite[3, Corollary 1]{On}) the sectional curvature of
$L \backslash G_{\theta}$  is non negative. Note that  $\frak l \subset 
[\frak g _{\theta} , \frak g _{\theta}] ^{\perp}$, hence the fibers of $\pi$ are totally 
geodesic since  $\nabla _T V =0$ for $T, V \in \frak l$ (see 
Proposition~\ref{flat}).

We summarize the above paragraphs as follows:

\begin{teo} \label{complete} Let $(G_{\theta},\{ J_{\alpha} \} , g)$ be a simply 
connected hyper-K\"ahler Lie group, so that
$G_{\theta}= \R^s \times (\R^k \ltimes_{\theta} \H^q)$, $k\leq q$, $s+k=4p$. Fix
 a connected closed 
abelian isotropic subgroup $L\subset \R^k$ acting on $G_{\theta}$ by the action $A$ as in 
\eqref{action} and  denote by $\pi : G_{\theta} \to L \backslash G_{\theta}$ 
 the associated Riemannian submersion. Then 
 \begin{enumerate}
\item the action $A$ of
$L$ on $G_{\theta}$ is free  and  preserves both, the
metric $g$ and the symplectic forms $\omega_{\alpha}$, $\alpha =
1,2,3$. The $L$-equivariant  moment map is
$\mu ^{\theta}= \mu ^{\theta} _1 i
+\mu ^{\theta}_2 j + \mu ^{\theta}_3 k,$ with $\mu ^{\theta}_{\alpha}$ given by
\[  \mu ^{\theta}_{\alpha}(X,W)( V) = \omega _{\alpha}(V,X) +
\frac 12\omega _{\alpha}(\rho _{\theta}(V) W, W),
\]
for any  $X \in \R^s \times \R^k, W \in \H^q, V \in \frak l$;
 \item $ L \backslash G_{\theta}$ has non negative sectional curvature. Moreover,  the fibers of  $\pi$ 
 are totally geodesic; 
\item for any $\xi \in  \text{Im}\, \mu ^{\theta} $,  
$ L \backslash {\mu ^{\theta}}^{-1}(\xi)$ is a closed embedded submanifold of   
$L \backslash G_{\theta}$ and  
the hyper-K\"ahler metric on the 
quotient $ L \backslash ({\mu ^{\theta}})^{-1}(\xi)$ is complete. \end{enumerate}
\end{teo}

\section{Examples}
In the next examples we show that it is possible to describe families
of known hyper-K\"ahler metrics \cite{GRG} in a unified way, by applying the
quotient construction to hyper-K\"ahler Lie groups $G$ with the action of a
  suitable abelian subgroup $L$ by left translations.

\subsection{Taub-Nut metric}

%It follows from Proposition \ref{flat} that
Let $\frak g_{\theta}$ be the  one parameter family of hyper-K\"ahler  
Lie algebras in dimension $8$ (see  
\eqref{8dim}) and $G_{\theta}$ the corresponding simply connected Lie groups.
  Let $\R$ be the subgroup of   $G_{\theta}= \H \ltimes _{\theta} \H$ given by
$(t,0)$, $t\in \R$, and let it act on $G_{\theta}$ by left translations, that is:
  \begin{eqnarray*}
\R \times G_{\theta} &\to &  \; \; G_{\theta} \\
(t ,(q,w))& \to &(t+ q, e^{i \theta t } w).
\end{eqnarray*}

Observe that Im$\,\H$ acts trivially on the second factor.
The corresponding hyper-K\"ahler moment map is
\begin{eqnarray*} \mu ^{\theta} & = &-\text{Im} (q )- \frac{\theta}2 \left( \text{Re} (iw
i\overline{w} )i +
\text{Re} (iw j\overline{w} )j + \text{Re} (iw k\overline{w} )k\right) \\
&=& -\text{Im} (q )+ \frac{\theta}2 \overline{w} i w.
\end{eqnarray*}
It can be checked that $\mu ^{\theta}$ is $L$-equivariant. 
The complete hyper-K\"ahler metric on $ \R \backslash (\mu ^{\theta})^{-1} (0)$ is the Taub-Nut
metric with parameter $\theta ^{-1} $ \cite{GRG}.

\subsection{Generalized Taubian-Calabi metric}

Let  $\theta =(\theta_1, \ldots , \theta _m)\in \R^m, \; G_{\theta}= \H  \ltimes _{\theta} \H^m$ and $\R =\{ (t,0)\, : \, t \in
\R \}$ acting on $G _{\theta}$ by left translations:
  \begin{eqnarray*}
\R \times G _{\theta} &\to &\;\; G _{\theta}\\
(t, (q,w_1, \ldots, w_m)) &\to &(t + q, e^{i\theta_1 t} w_1,
\ldots, e^{i \theta _m t} w_m).
\end{eqnarray*}
For $m = 1$ this is the Lie group considered in the
first example.
The corresponding hyper-K\"ahler moment map is
\begin{eqnarray*} \mu ^{\theta} & = & -\text{Im} (q )- \frac 12 \sum_{\beta =
1}^m \theta_{\beta} (  \text{Re} (iw_{\beta} i\overline{w_{\beta}} )i +
\text{Re} (iw_{\beta} j\overline{w_{\beta}} )j + \text{Re} (iw_{\beta} k\overline{w_{\beta}} )k) \\
{}&=& -\text{Im} (q )+ \frac 12 \sum_{{\beta} = 1}^m  \theta_{\beta}\overline{w_{\beta}} i w_{\beta}.
\end{eqnarray*}
When $\theta_{\beta}=1$, for each ${\beta}$, the complete hyper-K\"ahler metric on $\R \backslash (\mu ^{\theta})^{-1} (0)$
coincides with the Taubian-Calabi metric \cite{Ro,GRG}.

\subsection{Lee-Weinberg-Yi metric}

   Let $\theta \in \text{GL}(m, \R), \; G _{\theta}= \H^m \ltimes _{\theta} \H^m$ and  $$\R^m =\{ ((t_1, \dots , t_m), 0)
\, : \, t_i \in \R \}$$
acting on $G _{\theta}$ by left translations:
  \begin{eqnarray*}
\R^m \times G _{\theta} &\to &\;\; G _{\theta} \\
((t_1, \ldots, t_m), (q_1, \ldots, q_m, w_1, \ldots, w_m)) &\to &(t_1 + q_1
, \ldots,  t_m + q_m ,\\{}&{}& e^{i \langle \theta_1 ,T \rangle} w_1, \ldots, 
e^{i \langle \theta_m , T \rangle } w_m), 
\end{eqnarray*}
where $T =(t_1, \ldots , t_m), \; \theta_{\beta}$ are the rows of $\theta$ and 
$\langle \,\, , \, \rangle$ is the Euclidean inner product in $\R^m$. 
The corresponding hyper-K\"ahler moment map is
$$
   \mu ^{\theta}  =  \left(- \text{Im} (q_1 )+ \frac 12 \sum_{\beta=1}^m \theta^1_{\beta} \overline{w_{\beta}} i w_{\beta} ,
\dots ,  - \text{Im} (q_m )+ \frac 12 \sum_{\beta=1}^m \theta^m_{\beta}\overline{w_{\beta}} i w_{\beta} \right).
$$
The  complete hyper-K\"ahler metric on $\R \backslash \mu^{-1}_{\theta} (0)$ is the
Lee-Weinberg-Yi metric with \newline
 $(\lambda_b^a) = \theta ^{-1}$
\cite{LWY,Mu,GRG}.

\section{Topology of the quotient  and local description of the metric}

  Let $G_{\theta}$ be the simply connected hyper-K\"ahler Lie group  $\R^s
\times (\R ^k\ltimes _{\theta} \H^q)$  ($s+k=4p$, $k \leq q$), $\theta \in 
M(k,q;k)$ and 
${\theta}(\H^q) = \H^{q}$. Let $L \subset \R^k$ be a closed abelian  subgroup
  with Lie algebra $\frak l = \text{span}_{\R}\{ e_1 , \ldots , e_l \}$
such that $\frak l$ is isotropic with respect 
to $\omega_{\alpha}$, for each $\alpha$.  
 
Let $T^q$ be the maximal torus of Sp$(q)$ with Lie algebra $\frak t ^q$ (defined in \eqref{torus}) 
 whose elements are of the form:
\begin{equation} \label{T^q} B = \begin{pmatrix}
  B (\phi_1)& &  \\
 &\ddots &  \\
 & & B(\phi_q)  \end{pmatrix},
\end{equation}
where $\phi_{\beta}\in \R$ and $B(\phi_{\beta})$ is the following $4\times 4$ real matrix: 
$$ \begin{pmatrix}
  \cos (\phi_{\beta})&-\sin (\phi_{\beta}) &0&0  \\
\sin (\phi_{\beta}) &\cos (\phi_{\beta}) &0&0 \\
0&0&\cos(\phi_{\beta})&-\sin(\phi_{\beta})  \\
0&0& \sin (\phi_{\beta})& \cos (\phi_{\beta}) \end{pmatrix}.
$$
We have an action $\varphi$ of $T^q$  on $G_{\theta}$:
\begin{equation} \label{conj}
\begin{array} {l}
\varphi: T^q \times G_{\theta} \to  G_{\theta}\\
(g, (X, W)) \to \varphi(g,(X,W)) = (X, B W),
\end{array}
\end{equation}
where $BW$ stands for the product of the $4q \times 4q$ matrix $B$ given in \eqref{T^q} 
by the column vector $W \in \H^q \cong \R^{4q}$. 
 Note that the action $\varphi$ commutes with $A$ (see \eqref{action}) and both,  
$A$ and $\varphi$, preserve the metric and are tri-holomorphic. Therefore, $T^q$  also 
 acts on the hyper-K\"ahler quotient by tri-holomorphic isometries.  Moreover, the next
theorem shows that when $l=p=q$, so that the quotient has dimension $4q$, the 
$T^q$-action has a unique fixed point.

In the next theorem we give the explicit description of the moment map 
 and show 
that the quotient manifold is diffeomorphic to the Euclidean space. 

\begin{teo} \label{hyperkahlerquotient} Let $G_{\theta}= \R ^s \times (\R^k \ltimes _{\theta} \H^q)$ be a hyper-K\"ahler Lie group, $s+k=4p$, 
$\theta \in M(k,q;k)$,   $L$ the connected closed abelian isotropic subgroup
$L \subset \R^k$  defined above and $A$, $\varphi$ 
 as in \eqref{action},  \eqref{conj}. Then 
 \begin{enumerate}
\item the expression of the moment map is  
\begin{equation}\label{explmu} \mu ^{\theta} (X,W) = 
\left(-{\text{Im}}\, X_1 +\frac 12 \sum_{\beta=1}^q \theta _{\beta}^1\overline{W_{\beta}} \, i W_{\beta},
 \dots , -{\text{Im}}\, X_l +\frac 12 \sum_{\beta=1}^q \theta _{\beta}^l\overline{W_{\beta}} \, i W_{\beta}\right), 
\end{equation}
for $(X,W) \in G_{\theta}$; 
\item we have the following diffeomorphisms:
\[  L \backslash (\mu ^{\theta})^{-1} (0) \cong \R^{4p+4q -  4l}, \qquad   
L\backslash  G_{\theta} \cong  \R^{4p +  4q -  l}; \]
\item  the torus $T^q$ acts on 
 $L \backslash
(\mu ^{\theta})^{-1} (0)$ by tri-holomorphic 
  isometries. If $l=p=q$,  the action of $T^q$ on the 
$4q$ dimensional quotient has  a unique fixed point.

\end{enumerate}
\end{teo}

\begin{proof}

%, since  
%$L \backslash \mu^{-1} (0)$ is a closed submanifold of the
%non negatively curved space $L \backslash G$. 
%The fibers of the riemannian submersion $G
%\stackrel{\pi}{\longrightarrow} L \backslash G$ are totally
%geodesic since $\nabla_XY=0$ if $X,Y \in \frak l$ (see Proposition~\ref{flat}). 
% \footnote{    estudiar la inmersion, o por ejemplo decir que siempre el
%cociente hiperk es una subvariedad riemanniana de un fibrado asi y
%asa}

In order to prove the second assertion  we will find global coordinates 
on $L \backslash \mu_{\theta}^{-1} (0)$ and $L\backslash  G_{\theta} $. For $(X,W) \in G_{\theta}$, $(T,0) \in L$, set
\begin{eqnarray} \label{basis}
X &= &\sum_{\alpha=1}^p e_{\alpha} (x_{\alpha} +  b_{\alpha}i + 
s_{\alpha}j +  p_{\alpha} k), \qquad T  =\sum_{\alpha = 1}^l t_{\alpha}
e_{\alpha},\\
   W& =& \sum_{\alpha=1}^q f_{\alpha} (u_{\alpha} +  y_{\alpha}i + 
z_{\alpha}j +  w_{\alpha}k).
   \end{eqnarray}
It follows that $(x_{\alpha}, b_{\gamma}, s_{\gamma}, p_{\gamma},
u_{\beta}, y_{\beta}, z_{\beta}, w_{\beta})$, with $\alpha = l+1, \ldots,
p$, $\gamma = l , \ldots, p$, \newline
 $\beta = 1, \ldots, q$, are global
coordinates on $L \backslash G_{\theta}$ and therefore $L\backslash  G_{\theta}$ is diffeomorphic to $\R^{4p +  4q -  l}$. 
Using the fact that the hypercomplex structure  corresponds to
$$
J_1 = R_{- i}, \quad J_2= R_{- j}, \quad J_3= R_{- k}
$$
    and that the metric $g$ is such that the real basis $$\{
e_{\alpha},  e_{\alpha} i, e_{\alpha} j,  e_{\alpha} k, f_{\beta},
f_{\beta} i,  f_{\beta} j,  f_{\beta} k, \;\; 1 \leq \alpha \leq p, \; 1 \leq
\beta \leq q \}$$ is orthonormal, we get
the following expression of the moment maps ${\mu ^{\theta}}_{\gamma}$, $\gamma =
1, 2, 3$, in terms of the real coordinates on $\H^p$ and $\H^q$:
\begin{eqnarray*}
{\mu ^{\theta}}_1 (X, W) (T) &=&  -\sum_{\alpha = 1}^l b_{\alpha} t_{\alpha}  +
\frac{1} {2}  \sum_{\alpha = 1}^l t_{\alpha} \left(\sum_{\beta
= 1}^q \theta_{\beta}^{\alpha} ( u_{\beta}^2 + y_{\beta}^2 -
z_{\beta}^2 - w_{\beta}^2)\right)\\ {}&=& g\left(T, \sum_{\alpha = 1}^l
\left(-b_{\alpha}  + \frac{1} {2}  \sum_{\beta = 1}^q
\theta_{\beta}^{\alpha} ( u_{\beta}^2  + y_{\beta}^2 - z_{\beta}^2 -
w_{\beta}^2)\right) e_{\alpha}\right),\\
{\mu ^{\theta}}_2 (X, W) (T) &=&  -\sum_{\alpha = 1}^l s_{\alpha} t_{\alpha} +  
\sum_{\alpha = 1}^l t_{\alpha} \left(\sum_{\beta = 1}^q
\theta_{\beta}^{\alpha} ( -u_{\beta} w_{\beta} + z_{\beta}
y_{\beta})\right)\\ {}&=&  g\left(T, \sum_{\alpha = 1}^l \left(-s_{\alpha}  +
  \sum_{\beta = 1}^q \theta_{\beta}^{\alpha} (- u_{\beta} w_{\beta}
+ z_{\beta} y_{\beta})\right) e_{\alpha}\right),\\
{\mu ^{\theta}}_3 (X, W) (T) &=& - \sum_{\alpha = 1}^l p_{\alpha} t_{\alpha} +
  \sum_{\alpha = 1}^l t_{\alpha} \left(\sum_{\beta = 1}^q
\theta_{\beta}^{\alpha} ( u_{\beta} z_{\beta} + w_{\beta}
y_{\beta})\right)\\ {}&=&  g\left(T, \sum_{\alpha = 1}^l \left(-p_{\alpha}  +
  \sum_{\beta = 1}^q \theta_{\beta}^{\alpha} ( u_{\beta} z_{\beta}
+ w_{\beta} y_{\beta})\right) e_{\alpha}\right), 
\end{eqnarray*}
or, equivalently, \eqref{explmu} holds. 
On $(\mu ^{\theta})^{-1}(0)$ one has the following relations:
\begin{eqnarray*}
b_{\alpha} & = & \frac{1} {2} \sum_{\beta = 1}^q
\theta_{\beta}^{\alpha} ( u_{\beta}^2 + y_{\beta}^2 - z_{\beta}^2 -
w_{\beta}^2),\\ 
s_{\alpha} &=& -   \sum_{\beta = 1}^q
\theta_{\beta}^{\alpha} (u_{\beta} w_{\beta} - z_{\beta} y_{\beta}),\\
p_{\alpha} &= &  \sum_{\beta = 1}^q \theta_{\beta}^{\alpha}
(u_{\beta} z_{\beta}+ w_{\beta} y_{\beta}),
\end{eqnarray*}
for any $\alpha = 1 \ldots, l$.
Thus, one has that $(x_{\alpha}, b_{\gamma}, s_{\gamma}, p_{\gamma},
u_{\beta}, y_{\beta}, z_{\beta}, w_{\beta})$, with $\alpha = 1, \ldots,
p$, $\gamma = l + 1, \ldots, p$, $\beta = 1, \ldots, q$, are global
coordinates on $(\mu ^{\theta})^{-1} (0)$.

Since the action of $\R^l$ leaves
invariant $x_{\gamma}, b_{\gamma}, s_{\gamma}, p_{\gamma}$, ($
\gamma \geq  l + 1$) and rotates the coordinates $u_{\beta},
y_{\beta}, z_{\beta}, w_{\beta}$, $ \beta =1, \ldots, q$, one has that
$(x_{\gamma}, b_{\gamma},s_{\gamma}, p_{\gamma}, u_{\beta},
y_{\beta}, z_{\beta}, w_{\beta})$, with $\gamma = l + 1, \ldots, p$,
$\beta = 1, \ldots, q$, are global coordinates on $L \backslash
(\mu ^{\theta})^{-1} (0)
$. It follows that the quotient space is diffeomorphic to $\R^{4p
+ 4q - 4l}$.

3. It follows from \eqref{conj} that:
\[  \varphi(B) (X,W) = (X, B W ),\] 
$B\in T^q, (X,W) \in G_{\theta} . $   
Since the moment map for the   action $\varphi$ is $T^q$-equivariant,  
$\varphi$ preserves $(\mu ^{\theta})^{-1} (0)$. In particular, the hyper-K\"ahler 
quotient  admits a
tri-holomorphic action of the torus $T^q$.  

Assume next that $l= p=q$, hence $\frak l \oplus J_1 \frak l \oplus J_2 \frak l 
\oplus J_3 \frak l = \R ^{4p}$. Let $\pi$ be the natural projection from $(\mu ^{\theta})^{-1}(0)$ onto $L \backslash (\mu
^{\theta} )^{-1}(0)$.  Then $\pi (X,W)$, $(X,W) \in (\mu ^{\theta})^{-1}(0)$, is a fixed 
point for the action of $T^q$ if and only if \[  
(V,0)\cdot (X,W)\cdot (V,0)^{-1}\cdot(X,W)^{-1} \in L
\]
for every $V \in \R^q$. We will show that $\pi(X,W)=\pi(0,0)$, that is, $\pi(0,0)$ is the 
unique  fixed point.  
Using \eqref{conj} and \eqref{product}  we calculate \[
(V,0) \cdot (X,W)\cdot (V,0)^{-1}\cdot(X,W)^{-1} = (V,0)\cdot (-V,W-(-V)\cdot W) = (0, V\cdot W - W) \]
which belongs to $L$ if and only if $V\cdot W = W$ for every $V \in \R^q$, hence $W=0$. 
Since $(X,0) \in (\mu ^{\theta})^{-1}(0)$ it follows that $\omega _{\alpha}(V, X)=0$
for every $V \in \frak l$, $\alpha=1,2,3$, and the assumption on $l$ implies that $X\in \frak l$. Therefore  
$\pi(X,0)=\pi(0,0)$, as asserted. 
\end{proof}

Using the fact that the quotient admits a tri-holomorphic $T^q$-action, we can obtain the local expression of the
hyper-K\"ahler  metric  in terms of the structure constants of the Lie group
$G_{\theta}$.

 Observe that, if $l = p$, the quotient has dimension $4q$, thus by \cite{LR,PP} the induced hyper-K\"ahler metric can be
locally written as follows 
\begin{equation} \label{PPmetric}
\frac{1}{4} H_{\beta \gamma} d {\bf r}_{\beta} \cdot d {\bf r}_{\gamma} + \frac{1}{4} H^{\beta \gamma} (d \tau_{\beta} +
 {\xom}_{\beta \delta} \cdot d {\bf r}_{\delta}) (d \tau_{\gamma} +
 {\xom}_{\gamma \epsilon} \cdot d {\bf r}_{\epsilon}), 
\end{equation}
where $\beta ,\gamma =1, \ldots, q$, $(H)^{\beta \gamma}$ is the inverse of the matrix $(H)_{\beta \gamma}$. The Killing vector
fields $\frac{\partial} {\partial \tau_{\beta}}$  generate
the $T^q$-action, $\psi_{\beta},{\bf r}_{\beta}$ are
defined as in \eqref{Gibbonsflat} and we assume Einstein summation convention. If $l < p$, the quotient  splits as Riemannian product
of the flat Euclidean space $\R^{4p - 4l}$  by a
$4q$-dimensional hyper-K\"ahler manifold with a tri-holomorphic $T^q$-action.

\begin{teo}  The local expression of the hyper-K\"ahler metric on the  quotient $L \backslash (\mu
^{\theta} )^{-1}(0)$ is  
$h = h_0 + h_1$, where $h_0$ is the Euclidean metric on $\R^{4p - 4l}$ and $h_1$ is given by \eqref{PPmetric}, with
\begin{equation} \label{matrix}
H_{\beta \gamma} = (\tilde \theta  \tilde \theta^t)_{\beta \gamma} + \frac{1}{r_{\beta}} \, \delta_{\beta \gamma}.
\end{equation}
$\tilde \theta$ is the $q \times l$ matrix  obtained from $\theta$ by deleting the last $p - l$ columns, $\tilde
\theta^t$ is its transpose and $r_{\beta} = \vert {\bf r}_{\beta} \vert$.
\end{teo}

\begin{proof}
The action $A$ given by \eqref{action} of $L$ on $G_{\theta}$ in the coordinates $(x_{\alpha}, b_{\alpha} , s_{\alpha} ,
p_{\alpha},
\psi_{\beta}, {\bf r}_{\beta} )$ is 
$$
\begin{array} {l}
L \times G_{\theta} \to G_{\theta},\\
(T, (X_{\alpha}, \psi_{\beta},
{\bf r}_{\beta} )) \to (X_{\alpha} + t_{\alpha}, \psi_{\beta} + 2 \langle \theta_{\beta}, T \rangle, {\bf r}_{\beta}),
\end{array}
$$
with $\alpha =1, \ldots, p$, $t_{\alpha} = 0$ for $\alpha > l$, $\beta = 1, \ldots, q$ and $\psi_{\beta},{\bf r}_{\beta}$ are
defined as in \eqref{Gibbonsflat}.

The previous action leaves $$\tau_{\beta}= \psi_{\beta} - 2 \sum_{\alpha=1}^{p} \theta^{\alpha}_{\beta}x_{\alpha}, \quad
\beta = 1, \ldots, q,$$ invariant and $\frac{\partial} {\partial \tau_{\beta}}$ are Killing vector fields for the quotient
hyper-K\"ahler metric and generate the $T^q$-action induced by $\eqref{conj}$.

On $\mu^{-1} (0)$ one has
$$
{{\mbox {Im}} X_{\alpha}} = \frac{1} {2} \sum_{\beta = 1}^q \theta^{\alpha}_{\beta} {\bf r}_{\beta}, \quad \alpha = 1, \ldots,
l,
$$
so the metric on $\mu^{-1} (0)$ is 
$$
\sum_{\alpha =1}^p dx_{\alpha}^2 + \sum_{\alpha = l + 1}^p  d ({\mbox {Im}} X_{\alpha})^2 + \frac{1} {4} \sum_{\alpha = 1}^l
\left(\sum_{\beta = 1}^q
\theta^{\alpha}_{\beta} d {\bf r}_{\beta} \right)^2 + \frac{1} {4} \sum_{\beta = 1}^q \left(  \frac{1} {r_{\beta}} d {\bf
r}_{\beta}^2 + r_{\beta} (d \psi_{\beta} + {\xom}_{\beta}
\cdot d {\bf r}_{\beta})^2 \right).
$$
Projecting orthogonally to the Killing vector fields $\frac{\partial} {\partial x_{\alpha}}$, $\alpha = 1, \ldots, l$, one gets
that locally the metric on the quotient is given by $h = h_0 + h_1$, where
$$
h_0 = \sum_{\alpha = l + 1}^p  (dx_{\alpha}^2 + d ({\mbox {Im}} X_{\alpha})^2)
$$
and $h_1$ is given by \eqref{PPmetric} with the matrix $H$ as in \eqref{matrix}.

\end{proof}

\end{document}